\documentclass[11pt]{article}%
\usepackage{amsmath}
\usepackage{amsfonts}
\usepackage{amssymb}
\usepackage{graphicx}
\usepackage{theorem}
\usepackage{makeidx}
\usepackage{verbatim}
\usepackage[colorlinks=true]{hyperref}%
\providecommand{\U}[1]{\protect\rule{.1in}{.1in}}
\newtheorem{thm}{Theorem}[section]

\newtheorem{pr}[thm]{Proposition}
\newtheorem{df}[thm]{Definition}
\newtheorem{rmk}[thm]{Remark}
\newtheorem{cor}[thm]{Corollary}
{\theorembodyfont{\upshape}
\newtheorem{examp}[thm]{Example}
}
\numberwithin{equation}{section} 
\begin{document}

\title{ }

\begin{center}
\vspace*{1.5cm}

\textbf{CALCULUS OF TANGENT\ SETS AND
DERIVATIVES\ OF\ SET\ VALUED\ MAPS\ UNDER METRIC SUBREGULARITY\ CONDITIONS}

\vspace*{1cm}

M. DUREA

{\small {Faculty of Mathematics, "Al. I. Cuza" University,} }

{\small {Bd. Carol I, nr. 11, 700506 -- Ia\c{s}i, Romania,} }

{\small {e-mail: \texttt{durea@uaic.ro}}}

\bigskip

R. STRUGARIU

{\small Department of Mathematics, "Gh. Asachi" Technical University, }

{\small {Bd. Carol I, nr. 11, 700506 -- Ia\c{s}i, Romania,} }

{\small {e-mail: \texttt{rstrugariu@tuiasi.ro}}}
\end{center}

\bigskip

\bigskip

\noindent{\small {\textbf{Abstract: }}In this paper we intend to give some
calculus rules for tangent sets in the sense of Bouligand and Ursescu, as well
as for corresponding derivatives of set-valued maps. Both first and second
order objects are envisaged and the assumptions we impose in order to get the
calculus are in terms of metric subregularity of the assembly of the initial
data. This approach is different from those used in alternative recent papers
in literature and allows us to avoid compactness conditions. A special
attention is paid for the case of perturbation set-valued maps which appear
naturally in optimization problems.}

\bigskip


\bigskip


\bigskip

\section{Introduction}

The importance of the tangent sets in the study of various mathematical
problems, including optimization, viability theory, and control theory is well
known: Chapter 4 in the comprehensive monograph \cite{AubFra} clearly
emphasizes and illustrates this idea which is the basic layer of many papers
in literature. On the other hand, the most useful theoretical constructions of
tangent sets are those called in \cite{AubFra} the contingent and the adjacent
tangent sets which we call here Bouligand and respectively Ursescu tangent
sets after the names of the mathematicians who firstly introduced these concepts.

Once one has a concept of tangent set to an arbitrary set in a normed vector
space, then one can construct a corresponding derivative for set-valued maps,
and after that the question of calculus rules for these objects arises
naturally. Several recent papers in literature are devoted to the study of
such calculus for Bouligand and Ursescu first and second order tangent sets:
we cite here \cite{LeeHuy}, \cite{LiM}, \cite{Li} and the references therein.
The starting point of this paper is the remark that, in general, the quoted
papers employ quite strong conditions on the initial data in order to get the
desired calculus. One speaks here about conditions concerning several types of
generalized compactness for the graphs of underlying set-valued maps and such
requirements are known to be quite strong in infinite dimensional spaces.
Notice as well that the second-order conditions used in \cite{Li} are also of
the same nature. It is true that there is a difference between the conditions
used in \cite{LeeHuy} (working on finite dimensional vector spaces for first
order objects) and \cite{Li} (working on general Banach spaces for second
order objects) as the authors of the later paper emphasize, but, however, the
techniques and the arguments are not so different. Another remark is that in
all these papers the approach inherits the main line of arguments and
techniques from \cite[Chapter 3]{Luc}.

In contrast, we prefer here to follow the way open in \cite[Chapter 4]{AubFra}
allowing, as we shall see, to use metric (sub)regularity assumptions which are
more adequate for infinite dimensional setting. The role of metric regularity
in the regularity theory of constraint systems is nowadays well known (see the
monographs \cite{KlaKum}, \cite{RocWet}, \cite{Mor2006}, \cite{DontRock2009b}
for the main results in these topics and for detailed discussions) and this
fact give credit to our method. Other differences with respect to the previous
papers in literature are as follows: we impose conditions on the assembly of
the initial data rather than separate condition for each object and the
calculus rules we obtain refer, in general, not to equalities but to those
inclusions which are shown to be the right ones to use in getting necessary
optimality conditions in some vector optimization problems.

The paper is organized as follows. The second section introduces the main
notations, concepts and auxiliary results needed in the rest of the paper. The
third section concentrates on some general results concerning the calculus of
the first and second order tangent sets under metric subregularity
assumptions. Then a calculus rule for the derivatives of the sum between
set-valued maps is given and the proof is based on the reduction of this
situation to the general case of tangent sets previously obtained. An
application of this calculus rule to vector optimization problems and an
example underlying the differences with respect to the requirements used in
the other papers in literature are given. The fourth section employs a similar
technique to the case of generalized perturbation maps and several inclusion
concerning the coderivatives of these maps are presented. After every
important result of the paper we present sufficient conditions in terms of
Fr\'{e}chet normal cones to ensure the fulfillment of the metric subregularity
condition imposed in hypotheses. The paper ends with a short section which
displays the main conclusion of this work.

\section{Preliminaries and tools}

In the sequel, we suppose that the involved spaces are Banach, unless
otherwise stated. In this setting, $B(x,r)$ and $D(x,r)$ denote the open and
the closed ball with center $x$ and radius $r,$ respectively. If $x\in X$ and
$A\subset X,$ one defines the distance from $x$ to $A$ as $d(x,A):=\inf
\{\left\Vert x-a\right\Vert \mid a\in A\}.$ As usual, we use the convention
$d(x,\emptyset)=\infty.$ For a non-empty set $A\subset X$ we put
$\operatorname*{cl}A$ and $\operatorname*{int}A$ for its topological closure
and interior, respectively. On a product of normed vector spaces we consider
the sum norm and the corresponding topology. Usually, the zero element of $X$
is denoted by $0_{X}.$

Let $F:X\rightrightarrows Y$ be a multifunction. The domain and the graph of
$F$ are denoted respectively by $\operatorname*{Dom}F:=\{x\in X\mid
F(x)\neq\emptyset\}$ and $\operatorname*{Gr}F:=\{(x,y)\in X\times Y\mid y\in
F(x)\}.$ If $A\subset X$ then $F(A):=\bigcup_{x\in A}F(x).$ The inverse
set-valued map of $F$ is $F^{-1}:Y\rightrightarrows X$ given by $F^{-1}%
(y):=\{x\in X\mid y\in F(x)\}$.

\bigskip

One says that $F$ is open at linear rate $c>0$ around $(\overline{x}%
,\overline{y})\in\operatorname*{Gr}F$ if there exist two neighborhoods
$U\in\mathcal{V}(\overline{x}),$ $V\in\mathcal{V}(\overline{y})$ and a
positive number $\varepsilon>0$ such that, for every $(x,y)\in
\operatorname*{Gr}F\cap(U\times V)$ and every $\rho\in(0,\varepsilon),$%
\[
B(y,\rho c)\subset F(B(x,\rho)).
\]

It is well known that this property is equivalent to the metric regularity
property of $F$ around $(\overline{x},\overline{y})$ which requires to exist
$a>0$ and two neighborhoods $U\in\mathcal{V}(\overline{x}),$ $V\in
\mathcal{V}(\overline{y})$ such that for every $u\in U$ and every $v\in V$ to
have
\[
d(u,F^{-1}(v))\leq ad(v,F(u)).
\]

Recall as well that $F$ is said to have the Aubin property at a point
$(\overline{x},\overline{y})\in\operatorname*{Gr}F$ if there exist $L>0,$
$r>0$ such that for all $x^{\prime},x^{\prime\prime}\in D(\overline{x},r)$ we
have
\[
F(x^{\prime})\cap D(\overline{y},r)\subset F(x^{\prime\prime})+L\left\Vert
x^{\prime}-x^{\prime\prime}\right\Vert D(0,1).
\]

\bigskip

Denote by $X^{\ast}$ the topological dual of $X.$ As announced, we shall use
in the main sections of the paper the Fr\'{e}chet normal cones in order to
write sufficient conditions for the fulfillment of the metric regularity
assumptions we use in our results. Here are the main facts which allow us to
do this.

\begin{df}
Let $X$ be a normed vector space, $S$ be a non-empty subset of $X$ and let
$x\in S.$ The Fr\'{e}chet normal cone to $S$ at $x$ is%
\begin{equation}
\widehat{N}(S,x):=\left\{  x^{\ast}\in X^{\ast}\mid\underset{u\overset
{S}{\rightarrow}x}{\lim\sup}\frac{x^{\ast}(u-x)}{\left\Vert u-x\right\Vert
}\leq0\right\}  . \label{eps-no}%
\end{equation}

\end{df}

\begin{df}
Let $F:X\rightrightarrows Y$ be a set-valued map and $(\overline{x}%
,\overline{y})\in\operatorname*{Gr}F.$ Then the Fr\'{e}chet coderivative at
$(\overline{x},\overline{y})$ is the set-valued map $\widehat{D}^{\ast
}F(\overline{x},\overline{y}):Y^{\ast}\rightrightarrows X^{\ast}$ given by
\[
\widehat{D}^{\ast}F(\overline{x},\overline{y})(y^{\ast}):=\{x^{\ast}\in
X^{\ast}\mid(x^{\ast},-y^{\ast})\in\widehat{N}(\operatorname{Gr}%
F,(\overline{x},\overline{y}))\}.
\]

\end{df}

If $F=f$ is a single-valued function, then we write $\widehat{D}^{\ast
}f(\overline{x})$ for $\widehat{D}^{\ast}F(\overline{x},\overline{y}).$ If $f$
is Fr\'{e}chet differentiable at $\overline{x},$ then $\widehat{D}^{\ast
}f(\overline{x})(y^{\ast})=\{\nabla^{\ast}f(\overline{x})(y^{\ast})\}$ for
every $y^{\ast}\in Y^{\ast},$ where $\nabla$ denotes the Fr\'{e}chet
differential, while $\nabla^{\ast}f(\overline{x})$ stands for the $[\nabla
f(\overline{x})]^{\ast}$ (the adjoint operator of $\nabla f(\overline{x})$)$.$

We recall a well-known openness result for set-valued maps (see, e.g.,
\cite[Theorem 5.6]{MorSha} or \cite[Theorem 4.1]{Mor2006}).

\begin{thm}
\label{desch3}Let $X,Y$ be Asplund spaces, $F:X\rightrightarrows Y$ be a
set-valued map with closed graph and $(\overline{x},\overline{y}%
)\in\operatorname*{Gr}F.$ Then the following assertions are equivalent:

(i) There exist $r>0,s>0$ and $c>0$ such that for every $(x,y)\in
\operatorname*{Gr}F\cap\lbrack B(\overline{x},r)\times B(\overline{y},s)]$ and
every $y^{\ast}\in Y^{\ast},x^{\ast}\in\widehat{D}^{\ast}F(x,y)(y^{\ast}),$
\begin{equation}
c\left\Vert y^{\ast}\right\Vert \leq\left\Vert x^{\ast}\right\Vert .
\label{surj}%
\end{equation}

(ii) There exist $\alpha>0,\beta>0,c>0$ and $\varepsilon>0$ such that for
every $(x,y)\in\operatorname*{Gr}F\cap\lbrack B(\overline{x},\alpha)\times
B(\overline{y},\beta)],$ every $a\in(0,c)$ and every $\rho\in(0,\varepsilon],$%
\[
B(y,\rho a)\subset F(B(x,\rho)).
\]

\end{thm}

From an inspection of the proof of this basic result, one observes that the
implication $(i)\Rightarrow(ii)$ remains true even outside the class of
Asplund spaces (see, for more details, \cite{DurStr}):

\begin{itemize}
\item If $\operatorname{Gr}F$ is convex if is enough that $X,Y$ to be Banach spaces.

\item If one replace the Fr\'{e}chet constructions with some other generalized
differentiation objects having reasonable good behavior on appropriate classes
of spaces.
\end{itemize}

We look in this paper to the case of the restriction of a single-valued map
$f:X\longrightarrow Y$ to a nonempty closed set $M\subset X$. We consider then
the set-valued map $F_{f,M}:X\rightrightarrows Y$ given by
\[
F_{f,M}(x)=\left\{
\begin{array}
[c]{l}%
\{f(x)\},x\in M\\
\emptyset,x\notin M.
\end{array}
\right.
\]

In fact, one can write $F_{f,M}=f+\Delta_{M},$ where $\Delta_{M}$ is the
indicator mapping of $M,$
\[
\Delta_{M}:X\rightrightarrows Y,\Delta_{M}(x)=\left\{
\begin{array}
[c]{l}%
\{0\}\subset Y,x\in M\\
\emptyset,x\notin M.
\end{array}
\right.  .
\]
It is well known (and easy to see) that $\widehat{D}^{\ast}\Delta
_{M}(\overline{x},0)(y^{\ast})=\widehat{N}(M,\overline{x})$ for every
$\overline{x}\in M$ and $y^{\ast}\in Y^{\ast}$.

Of course, $\operatorname*{Gr}F_{f,M}=\operatorname*{Gr}f\cap(M\times Y)$ and
in order to write condition (\ref{surj}) for $F_{f,M}$ we need some calculus rules.

One such calculus is proved in \cite[Theorem 4.3]{DurStr2}. Let
$f:X\rightarrow Y$ be calm at $\overline{x},$ i.e. there exist $l>0$ and a
neighborhood $U$ of $\overline{x}$ such that%
\[
\left\Vert f(x)-f(\overline{x})\right\Vert \leq l\left\Vert x-\overline
{x}\right\Vert \text{ for every }x\in U.
\]
and $G:X\rightrightarrows Y$ be an arbitrary set-valued map such that
$(\overline{x},{\overline{y}})\in\operatorname*{Gr}G$. Then%
\begin{equation}
\hat{D}^{\ast}(G-f)(\overline{x},{\overline{y}-f(}\overline{x}))(y^{\ast
})\subset%
%TCIMACRO{\dbigcap \limits_{x^{\ast}\in\hat{D}^{\ast}f(\overline{x})(y^{\ast}%
%)}}%
%BeginExpansion
{\displaystyle\bigcap\limits_{x^{\ast}\in\hat{D}^{\ast}f(\overline{x}%
)(y^{\ast})}}
%EndExpansion
[\hat{D}^{\ast}G(\overline{x},{\overline{y}})(y^{\ast})-x^{\ast}]\ \text{for
all }y^{\ast}\in Y^{\ast}, \label{dif_Ff}%
\end{equation}
provided that $\hat{D}^{\ast}f(\overline{x})(y^{\ast})\not =\emptyset.$
Furthermore, inclusion (\ref{dif_Ff}) holds as equality if $f$ is Fr\'{e}chet
differentiable at $\overline{x}$. Summing up all the facts already presented,
for every $\overline{x}\in M,$ if $f$ is Fr\'{e}chet differentiable at
$\overline{x},$
\[
\hat{D}^{\ast}F_{f,M}(\overline{x},{f(}\overline{x}))(y^{\ast})=\widehat
{N}(M,\overline{x})+\nabla^{\ast}f(\overline{x})(y^{\ast}),\forall y^{\ast}\in
Y^{\ast}.
\]

This formula and Theorem \ref{desch3} allow us to give the following consequence.

\begin{cor}
\label{cor_mr}Let $X,Y$ be Asplund spaces, $f:X\longrightarrow Y$ be a
Fr\'{e}chet differentiable function, $M\subset X$ be a closed set and
$\overline{x}\in M.$ Suppose that the following assumption is satisfied: there
exist $c>0,$ $r>0$ such that for every $x\in M\cap B(\overline{x},r)$ and
every $y^{\ast}\in Y^{\ast},$ $x^{\ast}\in\widehat{N}(M,x)+\nabla^{\ast
}f(x)(y^{\ast}),$%
\[
c\left\Vert y^{\ast}\right\Vert \leq\left\Vert x^{\ast}\right\Vert .
\]

Then $F_{f,M}$ is metrically regular around $(\overline{x},f(\overline{x})),$
that is, there exist $s>0,$ $\mu>0$ s.t. for every $u\in B(\overline{x},s)\cap
M$ and $v\in B(f(\overline{x}),s)$
\[
d(u,f^{-1}(v)\cap M)\leq\mu\left\Vert v-f(u)\right\Vert .
\]

\end{cor}

The property displayed in the above conclusion is met in literature under the
name of metric regularity of $f$ around $(\overline{x},f(\overline{x}))$ with
respect to $M$, a terminology we use here as well.

\bigskip

However, we use as well in this paper a weaker condition (see, e.g.,
\cite{ArtMord2010}, \cite{ZheNg}): one says that a set-valued map
$F:X\rightrightarrows Y$ is metrically subregular at $(\overline{x}%
,\overline{y})\in\operatorname*{Gr}F$ if there exist $a>0$ and $U\in
\mathcal{V}(\overline{x}),$ such that for every $u\in U$
\[
d(u,F^{-1}(\overline{y}))\leq ad(\overline{y},F(u)).
\]
With the above notation, we say that the function $f$ is metrically subregular
at $(\overline{x},f(\overline{x}))$ with respect to $M$ $(\overline{x}\in M)$
if $F_{f,M}$ is metrically subregular at $(\overline{x},f(\overline{x})).$

More precisely, $f$ is metrically subregular at $(\overline{x},f(\overline
{x}))$ with respect to $M$ $(\overline{x}\in M)$ if there exist $s>0,$ $\mu>0$
s.t. for every $u\in B(\overline{x},s)\cap M$
\[
d(u,f^{-1}(f(\overline{x}))\cap M)\leq\mu\left\Vert f(\overline{x}%
)-f(u)\right\Vert .
\]

Of course, metric regularity around a point implies metric subregularity at
that point. For several aspects on the links and the differences between these
notions the reader is refered to \cite[p. 178]{Mor2006} and \cite[Section
3H]{DontRock2009b}. For sufficient conditions of subregularity in terms of
coderivative one can consult the recent paper \cite{ZheNg}, but in order to
avoid too many technicalities, we use in the sequel the specialization in
various particular cases of the condition in Corollary \ref{cor_mr} which
implies metric regularity and whence metric subregularity of the underlying functions.

\bigskip

We introduce some other main definitions we use in the sequel. The objects are
standard, but they are used under different names in literature (see, for
instance, \cite{AubFra}).

\begin{df}
\label{defcon}Let $D$ be a nonempty subset of $X$ and $\overline{x}\in X$.

(i) The first order Bouligand tangent cone to $D$ at $\overline{x}$ is the
set
\[
T_{B}(D,\overline{x})=\{u\in X\mid\exists(t_{n})\downarrow0,\exists
(u_{n})\rightarrow u,\forall n\in\mathbb{N},\overline{x}+t_{n}u_{n}\in D\}
\]
where $(t_{n})\downarrow0$ means $(t_{n})\subset(0,\infty)$ and $(t_{n}%
)\rightarrow0.$

(ii) If $x_{1}\in X$ the second order Bouligand tangent set to $D$ at
$(\overline{x},x_{1})$ is the set
\[
T_{B}^{2}(D,\overline{x},x_{1})=\{u\in X\mid\exists(t_{n})\downarrow
0,\exists(u_{n})\rightarrow u,\forall n\in\mathbb{N},\overline{x}+t_{n}%
x_{1}+t_{n}^{2}u_{n}\in D\}.
\]

(iii) The first order Ursescu tangent cone to $D$ at $\overline{x}$ is the
set
\[
T_{U}(D,\overline{x})=\{u\in X\mid\forall(t_{n})\downarrow0,\exists
(u_{n})\rightarrow u,\forall n\in\mathbb{N},\overline{x}+t_{n}u_{n}\in D\}.
\]

(iv) If $x_{1}\in X$\ the second order Ursescu tangent set to $D$\ at
$(\overline{x},x_{1})$\ is the set
\[
T_{U}^{2}(D,\overline{x},x_{1})=\{u\in X\mid\forall(t_{n})\downarrow
0,\exists(u_{n})\rightarrow u,\exists n_{0}\in\mathbb{N},\forall n\geq
n_{0},\overline{x}+t_{n}x_{1}+t_{n}^{2}u_{n}\in D\}.
\]

\end{df}

\begin{rmk}
\label{obscon} For both types of tangent sets ($\ast\in\{B,U\}$) one has that
$T_{\ast}^{2}(D,\overline{x},x_{1})=\emptyset$ if $x_{1}\notin T_{\ast
}(D,\overline{x})$ and, in general, $T_{\ast}^{2}(D,\overline{x},0)=T_{\ast
}(D,\overline{x})$. In the same notation, the sets $T_{\ast}(D,\overline{x})$
are closed cones (not necessarily convex) but, in general, $T_{\ast}%
^{2}(D,\overline{x},x_{1})$ are not cones. Moreover, $T_{\ast}(D,\overline
{x})=T_{\ast}(\operatorname{cl}D,\overline{x})$ and $T_{\ast}^{2}%
(D,\overline{x},x_{1})=T_{\ast}^{2}(\operatorname{cl}D,\overline{x},x_{1})$.
If $\overline{x}\in\operatorname{int}A,$ then, $T_{\ast}(D,\overline
{x})=T_{\ast}(D\cap A,\overline{x})$. It is also clear that $T_{B}%
(D,\overline{x})\neq\emptyset$ if and only if $\overline{x}\in
\operatorname*{cl}D$ and, obviously, $T_{U}(D,\overline{x})\subset
T_{B}(D,\overline{x}).$
\end{rmk}

The next well-known properties (see \cite[Tables 4.1 and 4.2]{AubFra}) are
important in the sequel.

\begin{pr}
\label{prop_con} Let $D\subset X,E\subset Y$ be closed sets, $x\in D,y\in
E,x_{1}\in X,y_{1}\in Y$ and $f:D\rightarrow Y$ be a Fr\'{e}chet
differentiable map. 

(i) Then
\begin{align*}
T_{U}(D,x)\times T_{B}(E,y)  & \subset T_{B}(D\times E,(x,y))\subset
T_{B}(D,x)\times T_{B}(E,y),\\
T_{U}(D,x)\times T_{U}(E,y)  & =T_{U}(D\times E,(x,y))
\end{align*}
and
\begin{align*}
T_{U}^{2}(D,x,x_{1})\times T_{B}^{2}(E,y,y_{1})  & \subset T_{B}^{2}(D\times
E,(x,y),(x_{1},y_{1}))\\
& \subset T_{B}^{2}(D,x,x_{1})\times T_{B}^{2}(E,y,y_{1}).
\end{align*}

(ii) If $x\in D\cap f^{-1}(E)$ then
\[
T_{B}(D\cap f^{-1}(E),x)\subset T_{B}(D,x)\cap\nabla f(x)^{-1}(T_{B}(E,f(x)))
\]
and
\[
T_{U}(D\cap f^{-1}(E),x)\subset T_{U}(D,x)\cap\nabla f(x)^{-1}(T_{U}(E,f(x))).
\]

(iii) If $M\subset X,\overline{x}\in M,x_{1}\in X$ and $A:X\rightarrow Y$ is a
bounded linear operator, then
\[
\operatorname*{cl}A(T_{B}(M,\overline{x}))\subset T_{B}(A(M),A(\overline{x})),
\]
and
\[
\operatorname*{cl}A(T_{B}^{2}(M,\overline{x},x_{1}))\subset T_{B}%
^{2}(A(M),A(\overline{x}),A(x_{1})).
\]

\end{pr}

\section{General calculus of tangent sets and derivatives}

Let us now give a calculus result for the above sets. The main line follows
the method implemented in \cite[Chapter 4]{AubFra}. We formulate now the
conditions which ensure the converse inclusions in Proposition \ref{prop_con}
(ii). These conditions are the more general ones met in the literature.

\begin{thm}
\label{th_teor}Let $X,Y$ be Banach spaces, $D\subset X,E\subset Y$ be closed
sets and $f:X\rightarrow Y$ be a continuously Fr\'{e}chet differentiable map
and $\overline{x}\in D\cap f^{-1}(E).$ Suppose that $g:X\times
Y\longrightarrow Y,$ $g(x,y):=f(x)-y$ is metrically subregular at
$(\overline{x},f(\overline{x}),0)$ with respect to $D\times E.$

\noindent Then%
\[
T_{B}(D,\overline{x})\cap\nabla f(\overline{x})^{-1}(T_{U}(E,f(\overline
{x})))\subset T_{B}(D\cap f^{-1}(E),\overline{x})
\]%
\[
T_{U}(D,\overline{x})\cap\nabla f(\overline{x})^{-1}(T_{B}(E,f(\overline
{x})))\subset T_{B}(D\cap f^{-1}(E),\overline{x})
\]%
\[
T_{U}(D,\overline{x})\cap\nabla f(\overline{x})^{-1}(T_{U}(E,f(\overline
{x})))=T_{U}(D\cap f^{-1}(E),\overline{x}).
\]
If moreover, $f$ is twice continuously differentiable, then for every
$x_{1}\in X$:%
\begin{align*}
&  T_{B}^{2}(D,\overline{x},x_{1})\cap\nabla f(\overline{x})^{-1}(T_{U}%
^{2}(E,f(\overline{x}),\nabla f(\overline{x})(x_{1}))-2^{-1}\nabla
^{2}f(\overline{x})(x_{1},x_{1}))\\
&  \subset T_{B}^{2}(D\cap f^{-1}(E),\overline{x},x_{1}),\\
&  T_{U}^{2}(D,\overline{x},x_{1})\cap\nabla f(\overline{x})^{-1}(T_{B}%
^{2}(E,f(\overline{x}),\nabla f(\overline{x})(x_{1}))-2^{-1}\nabla
^{2}f(\overline{x})(x_{1},x_{1}))\\
&  \subset T_{B}^{2}(D\cap f^{-1}(E),\overline{x},x_{1}),\\
&  T_{U}^{2}(D,\overline{x},x_{1})\cap\nabla f(\overline{x})^{-1}(T_{U}%
^{2}(E,f(\overline{x}),\nabla f(\overline{x})(x_{1}))-2^{-1}\nabla
^{2}f(\overline{x})(x_{1},x_{1}))\\
&  \subset T_{U}^{2}(D\cap f^{-1}(E),\overline{x},x_{1}).
\end{align*}

\end{thm}

\noindent\textit{Proof.} According to the metric subregularity assumption,
there exists $s>0,$ $\mu>0$ s.t. for every $(x,y)\in\lbrack B(\overline
{x},s)\times B(f(\overline{x}),s)]\cap(D\times E)$
\begin{equation}
d((x,y),g^{-1}(0)\cap(D\times E))\leq\mu\left\Vert f(x)+y\right\Vert .
\label{met_reg_t}%
\end{equation}
Take $u\in T_{B}(D,\overline{x})\cap\nabla f(\overline{x})^{-1}(T_{U}%
(E,f(\overline{x}))),$ i.e. $u\in T_{B}(D,\overline{x})$ and $\nabla
f(\overline{x})(u)\in T_{U}(E,f(\overline{x})).$ Then there exist
$(t_{n})\downarrow0,(u_{n})\rightarrow u,$ $v_{n}\rightarrow\nabla f(x)(u)$
with $\overline{x}+t_{n}u_{n}\in D$ and $f(\overline{x})+t_{n}v_{n}\in E$ for
all $n$ large enough. Then one can apply $(\ref{met_reg_t})$ for every pair
$(u,v)=(\overline{x}+t_{n}u_{n},f(\overline{x})+t_{n}v_{n}).$ Then for every
$n$ large enough, there exists $(p_{n},q_{n})\in D\times E$ with
$f(p_{n})=q_{n}\ $and
\[
\left\Vert (\overline{x}+t_{n}u_{n},f(\overline{x})+t_{n}v_{n})-(p_{n}%
,q_{n})\right\Vert <\mu\left\Vert f(\overline{x}+t_{n}u_{n})-f(\overline
{x})-t_{n}v_{n}\right\Vert +t_{n}^{2}.
\]
Then for every $n$ as above, $p_{n}\in D\cap f^{-1}(E)$ and
\[
\left\Vert \overline{x}+t_{n}u_{n}-p_{n}\right\Vert <\mu\left\Vert
f(\overline{x}+t_{n}u_{n})-f(\overline{x})-t_{n}v_{n}\right\Vert +t_{n}^{2}%
\]
whence%
\[
\left\Vert t_{n}^{-1}(p_{n}-\overline{x})-u_{n}\right\Vert <\mu\left\Vert
t_{n}^{-1}[f(\overline{x}+t_{n}u_{n})-f(\overline{x})]-v_{n}\right\Vert
+t_{n}.
\]
Since $t_{n}^{-1}[f(\overline{x}+t_{n}u_{n})-f(\overline{x})]\overset
{n\rightarrow\infty}{\longrightarrow}\nabla f(\overline{x})(u)$, we infer that
$u_{n}^{\prime}:=t_{n}^{-1}(p_{n}-\overline{x})\rightarrow u$ which allows us
to conclude the proof of the first inclusion of the theorem. Now, the other
two first-order relations are similar. Notice that for the equality in the
third relation one take into account Proposition \ref{prop_con} (ii).

The proof of the second-order relations is similar. Nevertheless, we
illustrate it with the first inclusion. Take $u\in T_{B}^{2}(D,\overline
{x},x_{1})\cap\nabla f(\overline{x})^{-1}(T_{U}^{2}(E,f(\overline{x}),\nabla
f(\overline{x})(x_{1}))-2^{-1}\nabla^{2}f(\overline{x})(x_{1},x_{1})),$ i.e.
$u\in T_{B}^{2}(D,\overline{x},x_{1})$ and $\nabla f(\overline{x}%
)(u)+2^{-1}\nabla^{2}f(\overline{x})(x_{1},x_{1})\in T_{U}^{2}(E,f(\overline
{x}),\nabla f(\overline{x})(x_{1})).$ Then there exist $(t_{n})\downarrow
0,(u_{n})\rightarrow u$ s.t. $\overline{x}+t_{n}x_{1}+t_{n}^{2}u_{n}\in D$ and
$v_{n}\rightarrow\nabla f(\overline{x})(u)+2^{-1}\nabla^{2}f(\overline
{x})(x_{1},x_{1})$ with $f(\overline{x})+t_{n}\nabla f(\overline{x}%
)(x_{1})+t_{n}^{2}v_{n}\in E$ for all $n$ large enough. Again, one applies
$(\ref{met_reg_t})$ for the pairs $(\overline{x}+t_{n}x_{1}+t_{n}^{2}%
u_{n},f(\overline{x})+t_{n}\nabla f(\overline{x})(x_{1})+t_{n}^{2}v_{n}).$ For
all $n$ large enough, there exists $(p_{n},q_{n})\in D\times E$ with
$f(p_{n})=q_{n}\ $and
\begin{align*}
&  \left\Vert (\overline{x}+t_{n}x_{1}+t_{n}^{2}u_{n},f(\overline{x}%
)+t_{n}\nabla f(\overline{x})(x_{1})+t_{n}^{2}v_{n})-(p_{n},q_{n})\right\Vert
\\
&  <\mu\left\Vert f(\overline{x}+t_{n}x_{1}+t_{n}^{2}u_{n})-f(\overline
{x})-t_{n}\nabla f(\overline{x})(x_{1})-t_{n}^{2}v_{n}\right\Vert +t_{n}^{3}.
\end{align*}
Then $p_{n}\in D\cap f^{-1}(E)\ $and
\begin{align*}
&  \left\Vert t_{n}^{-2}(p_{n}-\overline{x}-t_{n}x_{1})-u_{n}\right\Vert \\
&  <\mu\left\Vert t_{n}^{-2}[f(\overline{x}+t_{n}x_{1}+t_{n}^{2}%
u_{n})-f(\overline{x})-t_{n}\nabla f(\overline{x})(x_{1})]-v_{n}\right\Vert
+t_{n}.
\end{align*}
Since $t_{n}^{-2}[f(\overline{x}+t_{n}x_{1}+t_{n}^{2}u_{n})-f(\overline
{x})-t_{n}\nabla f(\overline{x})(x_{1})]\rightarrow\nabla f(\overline
{x})(u)+2^{-1}\nabla^{2}f(\overline{x})(x_{1},x_{1})$ we get that
$u_{n}^{\prime}$ $:=t_{n}^{-2}(p_{n}-\overline{x}-t_{n}x_{1})\rightarrow u$
and $\overline{x}+t_{n}x_{1}+t_{n}^{2}u_{n}^{\prime}=p_{n}\in D\cap
f^{-1}(E).$

Then the proof concludes here.\hfill$\square$

\begin{rmk}
Following Corollary \ref{cor_mr}, a sufficient condition for the metric
(sub)regularity assumption in Theorem \ref{th_teor} could be written down (on
Asplund spaces) as follows: there exist $c>0,$ $r>0$ such that for every
$(x,y)\in(D\times E)\cap\lbrack B(\overline{x},r)\times B(\overline{y},r)]$
and every $y^{\ast}\in Y^{\ast},$ $(u^{\ast},v^{\ast})\in\widehat
{N}(D,x)\times\widehat{N}(E,y)+(\nabla^{\ast}f(x)(y^{\ast}),-y^{\ast})$%
\[
c\left\Vert y^{\ast}\right\Vert \leq\left\Vert (u^{\ast},v^{\ast})\right\Vert
.
\]

\end{rmk}

\begin{df}
\label{defder}Let $(\overline{x},\overline{y})\in\operatorname*{Gr}F.$ The
first order Bouligand derivative of $F$ at $(\overline{x},\overline{y})$ is
the set valued map $D_{B}F(\overline{x},\overline{y})$ from $X$ into $Y$
defined by
\[
\operatorname*{Gr}D_{B}F(\overline{x},\overline{y})=T_{B}(\operatorname*{Gr}%
F,(\overline{x},\overline{y})),
\]
and if $(x_{1},y_{1})\in X\times Y,$\ the second order Bouligand derivative of
$F$\ at $(\overline{x},\overline{y})$\ with respect to $(x_{1},y_{1})$\ is the
set valued map $D_{B}^{2}F((\overline{x},\overline{y}),(x_{1},y_{1}))$\ from
$X$\ into $Y$\ defined by
\[
\operatorname*{Gr}D_{B}^{2}F((\overline{x},\overline{y}),(x_{1},y_{1}%
))=T_{B}^{2}(\operatorname*{Gr}F,(\overline{x},\overline{y}),(x_{1},y_{1})).
\]

\end{df}

Now the first and second order Ursescu derivative has similar definition.

The first part of the following definition was introduced by J.-P. Penot
\cite{Pen}.

\begin{df}
\label{defdini}Let $(\overline{x},\overline{y})\in\operatorname*{Gr}F.$ (i)
The Dini lower derivative of $F$ at $(\overline{x},\overline{y})$ is the
multifunction $D_{D}F(\overline{x},\overline{y})$ from $X$ into $Y$ given, for
every $u\in X,$ by
\begin{align*}
D_{D}F(\overline{x},\overline{y})(u)  &  =\{v\in Y\mid\forall(t_{n}%
)\downarrow0,\forall(u_{n})\rightarrow u,\exists(v_{n})\rightarrow v,\exists
n_{0}\in\mathbb{N}\text{,}\\
\forall n  &  \geq n_{0},\overline{{y}}+t_{n}v_{n}\in F(\overline{x}%
+t_{n}u_{n})\}.
\end{align*}

(ii) \textit{If }$(x_{1},y_{1})\in X\times Y,$\textit{\ the second order Dini
lower derivative of }$F$ \textit{at }$(\overline{x},\overline{y}%
)$\textit{\ with respect to }$(x_{1},y_{1})$\textit{\ }is the multifunction
$D_{D}^{2}F(\overline{x},\overline{y})(x_{1},y_{1})$ from $X$ into $Y$ given,
for every $u\in X,$ by
\begin{align*}
D_{D}^{2}F((\overline{x},\overline{y}),(x_{1},y_{1}))(u) &  =\{v\in
Y\mid\forall(t_{n})\downarrow0,\forall(u_{n})\rightarrow u,\exists
(v_{n})\rightarrow v,\exists n_{0}\in\mathbb{N}\text{,}\\
\forall n &  \geq n_{0},\overline{y}+t_{n}y_{1}+t_{n}^{2}v_{n}\in
F(\overline{x}+t_{n}x_{1}+t_{n}^{2}u_{n})\}.
\end{align*}

\end{df}

\begin{rmk}
\label{obsincder}Obviously, the next inclusions are true for all
$(\overline{x},\overline{y})\in\operatorname*{Gr}F,$ $(x_{1},y_{1})\in X\times
Y$ and $u\in X$:
\begin{align*}
D_{D}F(\overline{x},\overline{y})(u) &  \subset D_{U}F(\overline{x}%
,\overline{y})(u)\subset D_{B}F(\overline{x},\overline{y})(u),\\
D_{D}^{2}F((\overline{x},\overline{y}),(x_{1},y_{1}))(u) &  \subset D_{U}%
^{2}F((\overline{x},\overline{y}),(x_{1},y_{1}))(u)\subset D_{B}%
^{2}F((\overline{x},\overline{y}),(x_{1},y_{1}))(u)
\end{align*}

\end{rmk}

One says that a set $A$ is derivable at a point $\overline{x}\in A$ if
$T_{B}(A,\overline{x})=T_{U}(A,\overline{x})$ (see \cite{AubFra}). Similarly,
one says that a set-valued map that $F$ is derivable (terminology of
\cite{AubFra}) or proto-differentiable (terminology of \cite{Roc}) at
$\overline{x}$ relative to $\overline{y}\in F(\overline{x})$ if its graph is
derivable at $(\overline{x},\overline{y}),$ i.e. $D_{U}F(\overline
{x},\overline{y})=D_{B}F(\overline{x},\overline{y}).$ One says that $F$ is
semi-differentiable at $\overline{x}$ relative to $\overline{y}\in
F(\overline{x})$ if $D_{D}F(\overline{x},\overline{y})=D_{B}F(\overline
{x},\overline{y}).$ It is clear that semi-differentiability implies
proto-differentiability. The next result is easy to prove.

\begin{pr}
\label{lpL}Suppose that $F$ has the Aubin property around $(\overline
{x},\overline{y})\in\operatorname*{Gr}F$ and $u\in X$\textit{. }Then:
\begin{align*}
D_{D}F(\overline{x},\overline{y})(u)=D_{U}F(\overline{x},\overline{y})(u)=\{
&  v\in Y\mid\forall(t_{n})\downarrow0,\exists(v_{n})\rightarrow v,\forall
n\in\mathbb{N},\\
&  \overline{y}+t_{n}v_{n}\in F(\overline{x}+t_{n}u)\}.
\end{align*}
Hence, if the set-valued map $F:X\rightrightarrows Y$ has the Aubin property
around $(\overline{x},\overline{y})\in\operatorname*{Gr}F$ and it is
proto-differentiable at $\overline{x}$ relative to $\overline{y}$, then it is
semi-differentiable at $\overline{x}$ relative to $\overline{y}.$
\end{pr}

However, as shown in \cite{Dur}, semi-differentiability is a quite strong
assumption which ensures some simple results concerning the derivative
calculus of the sum of set-valued maps. More explicitly, if $F_{1}%
,F_{2}:X\rightrightarrows Y$ are set-valued maps, for the sum $F_{1}%
+F_{2}:X\rightrightarrows Y$ given by
\begin{align*}
(F_{1}+F_{2})(x)  &  =F_{1}(x)+F_{2}(x)\\
&  =\{y\in Y\mid\exists y_{1}\in F_{1}(x),\exists y_{2}\in F_{2}%
(x),y=y_{1}+y_{2}\},
\end{align*}
one can easily prove (see \cite{Dur}) that for $(\overline{x},y_{1}%
)\in\operatorname*{Gr}F_{1},$ $(\overline{x},y_{2})\in\operatorname*{Gr}F_{2}$
if either $F_{1}$ is semi-differentiable at $\overline{x}$ relative to
$\overline{y}_{1}$ or $F_{2}$ is semi-differentiable at $\overline{x}$
relative to $\overline{y}_{2},$ then, for every $u\in X,$
\[
D_{B}F_{1}(\overline{x},\overline{y}_{1})(u)+D_{B}F_{2}(\overline{x}%
,y_{2})(u)\subset D_{B}(F_{1}+F_{2})(\overline{x},\overline{y}_{1}%
+\overline{y}_{2})(u)
\]
and%
\[
D_{U}F_{1}(\overline{x},\overline{y}_{1})(u)+D_{U}F_{2}(\overline{x}%
,\overline{y}_{2})(u)\subset D_{U}(F_{1}+F_{2})(\overline{x},\overline{y}%
_{1}+\overline{y}_{2})(u).
\]

\medskip

Despite the fact that it is a quite heavy assumption, this concept of
semi-differentiability is employed (in conjunction with some compactness
requirements) as the basic ingredient in several recent papers concerning the
calculus rules for both first and second order derivatives: see, for instance
\cite{Li}, \cite{LiM} and the references therein. We prefer here to avoid such
assumptions and to work with the far more natural hypotesis of proto-differentiability.

\bigskip

Before passing to the main results, we give some similar definitions for
second-order objects: one says that $A$ is second-order derivable at
$\overline{x}$ in the direction $x_{1}$ if $T_{B}^{2}(A,\overline{x}%
,x_{1})=T_{U}^{2}(A,\overline{x},x_{1}),$ and one says that $F$ is
second-order proto-differentiable at $\overline{x}$ relative to $\overline
{y}\in F(\overline{x})$ in the direction $(x_{1},y_{1})$ if $D_{B}%
^{2}F((\overline{x},\overline{y}),(x_{1},y_{1}))=D_{U}^{2}F((\overline
{x},\overline{y}),(x_{1},y_{1})).$

\medskip

Let $F:X\rightrightarrows Y$ be a set-valued map and $f:X\rightarrow X$ be a
single-valued map. Then we denote by $F\circ f$ the set-valued map from $X$ to
$Y$ given by $\left(  F\circ f\right)  (x)=F(f(x))$ for every $x\in X.$ We are
ready to present the first main result of the paper. The proof is based on a
transformation method which allows us to reduce the calculus of the involved
objects to the pattern of Theorem \ref{th_teor}.

\begin{thm}
\label{T1}Let $F_{1},F_{2}:X\rightrightarrows Y$ be set-valued maps with
closed graph, $f:X\rightarrow X$ be a continuously differentiable
single-valued map and $(\overline{x},\overline{y}_{1})\in\operatorname*{Gr}%
F_{1},$ $(\overline{x},\overline{y}_{2})\in\operatorname*{Gr}(F_{2}\circ f)$.
Suppose that the function $g:(X\times Y)^{2}\longrightarrow X,$ $g(\alpha
,\beta,\gamma,\delta)=f(\alpha)-\gamma$ is metrically subregular at
$(\overline{x},\overline{y}_{1},f(\overline{x}),\overline{y}_{2},0_{X})$ with
respect to $\operatorname*{Gr}F_{1}\times\operatorname*{Gr}F_{2}$.

(i) If either $F_{1}$ is proto-differentiable at $\overline{x}$ relative to
$\overline{y}_{1}$ or $F_{2}$ is proto-differentiable at $f(\overline{x})$
relative to $\overline{y}_{2},$ then, for every $u\in X,$
\[
D_{B}F_{1}(\overline{x},\overline{y}_{1})(u)+D_{B}F_{2}(f(\overline
{x}),\overline{y}_{2})(\nabla f(\overline{x})(u))\subset D_{B}(F_{1}%
+F_{2}\circ f)(\overline{x},\overline{y}_{1}+\overline{y}_{2})(u).
\]

(ii) Let $x\in X$ and $y_{1},y_{2}\in Y.$ If $f$ is linear and either $F_{1}$
is second-order proto-differentiable at $\overline{x}$ relative to
$\overline{y}_{1}$ in the direction $(x,y_{1})$ or $F_{2}$ is second-order
proto-differentiable at $\overline{x}$ relative to $\overline{y}_{2}$ in the
direction $(f(x),y_{2}),$ then, for every $u\in X,$%
\begin{align*}
&  D_{B}^{2}F_{1}((\overline{x},\overline{y}_{1}),(x,y_{1}))(u)+D_{B}^{2}%
F_{2}((f(\overline{x}),\overline{y}_{2}),(f(x),y_{2}))(\nabla f(\overline
{x})(u))\\
&  \subset D_{B}^{2}(F_{1}+F_{2}\circ f)((\overline{x},\overline{y}%
_{1}+\overline{y}_{2}),(f(x),y_{1}+y_{2}))(u).
\end{align*}

\end{thm}

\noindent\textbf{Proof.} (i) Let us consider the following auxiliary
functions:
\[
\varphi:\left(  X\times Y\right)  ^{2}\rightarrow X\times Y,\quad
\varphi(\alpha,\beta,\gamma,\delta)=(\alpha,\beta+\delta)
\]
and
\[
\psi:\left(  X\times Y\right)  ^{2}\rightarrow X,\quad\psi(\alpha,\beta
,\gamma,\delta)=f(\alpha)-\gamma.
\]
It is not difficult to see that $\operatorname*{Gr}(F_{1}+F_{2}\circ
f)=\varphi\left(  (\operatorname*{Gr}F_{1}\times\operatorname*{Gr}F_{2}%
)\cap\psi^{-1}(0_{X})\right)  \ $because
\begin{align*}
\varphi\left(  (\operatorname*{Gr}F_{1}\times\operatorname*{Gr}F_{2})\cap
\psi^{-1}(0_{X})\right)   &  =\varphi\left(  (\operatorname*{Gr}F_{1}%
\times\operatorname*{Gr}F_{2})\cap\{(x,y,z,t)\in\left(  X\times Y\right)
^{2}\mid f(x)=z\}\right) \\
&  =\varphi(\{(x,y,f(x),t)\mid y\in F_{1}(x),t\in F_{2}(f(x))\})\\
&  =\{(x,y+t)\mid x\in X,y+t\in(F_{1}+F_{2}\circ f)(x)\}\\
&  =\operatorname*{Gr}(F_{1}+F_{2}\circ f).
\end{align*}
The linearity of $\varphi$ and Proposition \ref{prop_con} $(iii)$, ensures:
\begin{align*}
\operatorname*{Gr}D_{B}(F_{1}+F_{2}\circ f)(\overline{x},\overline{y}%
_{1}+\overline{y}_{2})  &  =T_{B}(\operatorname*{Gr}(F_{1}+F_{2}\circ
f),(\overline{x},\overline{y}_{1}+\overline{y}_{2}))\\
&  =T_{B}(\varphi\left(  (\operatorname*{Gr}F_{1}\times\operatorname*{Gr}%
F_{2})\cap\psi^{-1}(0_{X})\right)  ,\varphi(\overline{x},\overline{y}%
_{1},f(\overline{x}),\overline{y}_{2}))\\
&  \supset\operatorname*{cl}\left(  \varphi(T_{B}((\operatorname*{Gr}%
F_{1}\times\operatorname*{Gr}F_{2})\cap\psi^{-1}(0_{X}),(\overline
{x},\overline{y}_{1},f(\overline{x}),\overline{y}_{2})))\right)  .
\end{align*}
Our hypotheses allow us to apply Theorem \ref{th_teor} for
$D:=\operatorname*{Gr}F_{1}\times\operatorname*{Gr}F_{2},$ $E:=\{0_{X}\}$ and
$f:=\psi.$ Because of the very particular form of $E,$ one observes that the
assumption on $g$ ensures that $h:(X\times Y\times X\times Y)\times
X\longrightarrow X,$ $h(\alpha,\beta,\gamma,\delta,\varepsilon)=f(\alpha
)-\gamma-\varepsilon$ is submetrically regular at $(\overline{x},\overline
{y}_{1},f(\overline{x}),\overline{y}_{2},0_{X},0_{X})$ with respect to
$(\operatorname*{Gr}F_{1}\times\operatorname*{Gr}F_{2})\times\{0_{X}\}$, so we
can indeed specialize Theorem \ref{th_teor} to the case described above.
Taking into account the equality $T_{U}(\{0_{X}\},0_{X})=T_{B}(\{0_{X}%
\},0_{X})=\{0_{X}\},$ we successively have
\begin{align*}
T_{B}((\operatorname*{Gr}F_{1}\times\operatorname*{Gr}F_{2})\cap\psi
^{-1}(0_{X}),(\overline{x},\overline{y}_{1},f(\overline{x}),\overline{y}%
_{2}))  &  \supset T_{B}(\operatorname*{Gr}F_{1}\times\operatorname*{Gr}%
F_{2},(\overline{x},\overline{y}_{1},f(\overline{x}),\overline{y}_{2}))\\
&  \cap\nabla\psi(\overline{x},\overline{y}_{1},f(\overline{x}),\overline
{y}_{2})^{-1}(T_{U}(\{0_{X}\},0_{X}))\\
&  =T_{B}(\operatorname*{Gr}F_{1}\times\operatorname*{Gr}F_{2},(\overline
{x},\overline{y}_{1},f(\overline{x}),\overline{y}_{2}))\\
&  \cap\nabla\psi(\overline{x},\overline{y}_{1},f(\overline{x}),\overline
{y}_{2})^{-1}(0_{X}).
\end{align*}
From Proposition \ref{prop_con} $(i)$ and the proto-differentiability
assumption, we get the following chain of inclusions:%
\begin{align*}
&  \operatorname*{Gr}D_{B}(F_{1}+F_{2}\circ f)(\overline{x},\overline{y}%
_{1}+\overline{y}_{2})\\
&  \supset\varphi(T_{B}(\operatorname*{Gr}F_{1}\times\operatorname*{Gr}%
F_{2},(\overline{x},\overline{y}_{1},f(\overline{x}),\overline{y}_{2}%
))\cap\nabla\psi(\overline{x},\overline{y}_{1},f(\overline{x}),\overline
{y}_{2})^{-1}(0_{X}))\\
&  =\varphi\left(  T_{B}(\operatorname*{Gr}F_{1},(\overline{x},\overline
{y}_{1}))\times T_{B}(\operatorname*{Gr}F_{2},(f(\overline{x}),\overline
{y}_{2})\right)  \cap\nabla\psi(\overline{x},\overline{y}_{1},f(\overline
{x}),\overline{y}_{2})^{-1}(0_{X}))\\
&  =\varphi(\{(u,v,w,p)\in(X\times Y)^{2}\mid(u,v)\in\operatorname*{Gr}%
D_{B}F_{1}(\overline{x},\overline{y}_{1}),\\
&  \hspace{1.9cm}(w,p)\in\operatorname*{Gr}D_{B}F_{2}(f(\overline
{x}),\overline{y}_{2}),\nabla f(\overline{x})(u)=w\})\\
&  =\{(u,v+p)\mid(u,v)\in\operatorname*{Gr}D_{B}F_{1}(\overline{x}%
,\overline{y}_{1}),(w,p)\in\operatorname*{Gr}D_{B}F_{2}(f(\overline
{x}),\overline{y}_{2}),\nabla f(\overline{x})(u)=w\}\\
&  =\operatorname*{Gr}(D_{B}F_{1}(\overline{x},\overline{y}_{1})(\cdot
)+D_{B}F_{2}(f(\overline{x}),\overline{y}_{2})(\nabla f(\overline{x}%
)(\cdot))).
\end{align*}

The proof of the first-order calculus rule is complete.

(ii) For the second part,
\begin{align*}
&  \operatorname*{Gr}D_{B}^{2}(F_{1}+F_{2}\circ f)((\overline{x},\overline
{y}_{1}+\overline{y}_{2}),(x,y_{1}+y_{2}))\\
&  =T_{B}^{2}(\operatorname*{Gr}(F_{1}+F_{2}\circ f),(\overline{x}%
,\overline{y}_{1}+\overline{y}_{2}),(x,y_{1}+y_{2}))\\
&  =T_{B}^{2}(\varphi\left(  (\operatorname*{Gr}F_{1}\times\operatorname*{Gr}%
F_{2})\cap\psi^{-1}(0)\right)  ,\varphi(\overline{x},\overline{y}%
_{1},f(\overline{x}),\overline{y}_{2}),\varphi(x,y_{1},f(x),y_{2}))\\
&  \supset\varphi\left(  T_{B}^{2}((\operatorname*{Gr}F_{1}\times
\operatorname*{Gr}F_{2})\cap\psi^{-1}(0),(\overline{x},\overline{y}%
_{1},f(\overline{x}),\overline{y}_{2}),(x,y_{1},f(x),y_{2}))\right)  .
\end{align*}

Now, we apply Theorem \ref{th_teor} for the same data as before. Taking into
account the equality $T_{U}^{2}(\{0_{X}\},0_{X},0_{X})=T_{B}^{2}%
(\{0_{X}\},0_{X},0_{X}))=\{0_{X}\}$ and because $\nabla f(\overline
{x})(x)=f(x)$ (since $f$ is linear), we successively have
\begin{align*}
&  T_{B}^{2}(\operatorname*{Gr}F_{1}\times\operatorname*{Gr}F_{2})\cap
\psi^{-1}(0_{X}),(\overline{x},\overline{y}_{1},f(\overline{x}),\overline
{y}_{2}),(x,y_{1},f(x),y_{2}))\\
&  \supset T_{B}^{2}(\operatorname*{Gr}F_{1}\times\operatorname*{Gr}%
F_{2},(\overline{x},\overline{y}_{1},f(\overline{x}),\overline{y}%
_{2}),(x,y_{1},f(x),y_{2}))\\
&  \cap\nabla\psi(\overline{x},\overline{y}_{1},f(\overline{x}),\overline
{y}_{2})^{-1}(0_{X}).
\end{align*}

Therefore, from the second-order proto-differentiability condition,
\begin{align*}
&  \operatorname*{Gr}D_{B}^{2}(F_{1}+F_{2}\circ f)((\overline{x},\overline
{y}_{1}+\overline{y}_{2}),(x,y_{1}+y_{2}))\\
&  \supset\varphi\left(  T_{B}^{2}(\operatorname*{Gr}F_{1}\times
\operatorname*{Gr}F_{2},(\overline{x},\overline{y}_{1},f(\overline
{x}),\overline{y}_{2}),(x,y_{1},f(x),y_{2}))\cap\nabla\psi(\overline
{x},\overline{y}_{1},f(\overline{x}),\overline{y}_{2})^{-1}(0_{X})\right)  \\
&  =\varphi\left(  T_{B}^{2}(\operatorname*{Gr}F_{1},(\overline{x}%
,\overline{y}_{1}),(x,y_{1}))\times T_{B}^{2}(\operatorname*{Gr}%
F_{2},(f(\overline{x}),\overline{y}_{2}),(f(x),y_{2})\right)  \cap\nabla
\psi(\overline{x},\overline{y}_{1},f(\overline{x}),\overline{y}_{2}%
)^{-1}(0_{X}))\\
&  =\varphi(\{(u,v,w,p)\in(X\times Y)^{2}\mid(u,v)\in\operatorname*{Gr}%
D_{B}^{2}F_{1}((\overline{x},\overline{y}_{1}),(x,y_{1})),\\
&  \hspace{1.9cm}(w,p)\in\operatorname*{Gr}D_{B}^{2}F_{2}((f(\overline
{x}),\overline{y}_{2})(f(x),y_{2})),\nabla f(\overline{x})(u)=w\}\\
&  =\{(u,v+p)\mid(u,v)\in\operatorname*{Gr}D_{B}^{2}F_{1}((\overline
{x},\overline{y}_{1}),(x,y_{1})),\\
&  \hspace{1.9cm}(w,p)\in\operatorname*{Gr}D_{B}^{2}F_{2}((f(\overline
{x}),\overline{y}_{2}),(f(x),y_{2})),\nabla f(\overline{x})(u)=w\}\\
&  =\operatorname*{Gr}(D_{B}^{2}F_{1}((\overline{x},\overline{y}_{1}%
),(x,y_{1}))(\cdot)+D_{B}^{2}F_{2}((f(\overline{x}),\overline{y}%
_{2}),(f(x),y_{2}))(\nabla f(\overline{x})(\cdot))).
\end{align*}

Consequently, the final formula is proved.$\hfill\square$

\begin{rmk}
Use again Corollary \ref{cor_mr} in order to get a sufficient condition for
the metric (sub)regularity assumption in Theorem \ref{T1} (on Asplund spaces):
there exist $c>0,$ $r>0$ such that for every $(x_{1},y_{1},x_{2},y_{2}%
)\in(\operatorname*{Gr}F_{1}\times\operatorname*{Gr}F_{2})\cap\lbrack
B(\overline{x},r)\times B(\overline{y},r)\times B(f(\overline{x}),r)\times
B(\overline{y},r)]$ and every $x^{\ast}\in X^{\ast},$ $(u_{1}^{\ast}%
,v_{1}^{\ast},u_{2}^{\ast},v_{2}^{\ast})\in\widehat{N}(\operatorname*{Gr}%
F_{1},(x_{1},y_{1}))\times\widehat{N}(\operatorname*{Gr}F_{2},(x_{2}%
,y_{2}))+(\nabla^{\ast}f(x_{1})(x^{\ast}),0_{Y^{\ast}},-x^{\ast},0_{Y^{\ast}%
})$%
\[
c\left\Vert x^{\ast}\right\Vert \leq\left\Vert (u_{1}^{\ast},v_{1}^{\ast
},u_{2}^{\ast},v_{2}^{\ast})\right\Vert .
\]

\end{rmk}

\bigskip

As one can see, we do not obtain equalities in Theorem \ref{T1}, but we do not
use any semi-differentiability or other compactness-like assumptions (compare
with \cite[Proposition 2]{LiM}). Nevertheless, we would like to emphasize that
the inclusions in Theorem \ref{T1} are exactly those needed in order to deduce
optimality conditions for some special types of vector optimization problems
in terms of the initial data. For instance, take $F_{1},F_{2}%
:X\rightrightarrows Y$ set-valued maps and $C\subset Y$ a closed convex cone
with nonempty interior which, as usual, gives a partial order relation on $Y.$
Consider the general vector optimization problem
\[
(P)\hspace{0.56in}\min(F_{1}+F_{2})\text{ s.t. }x\in X.
\]
{A point} $(\overline{x},\overline{y})\in\operatorname{Gr}(F_{1}+F_{2})$ {is
called a weak Pareto solution of }$(P)$ {if} { }%
\[
((F_{1}+F_{2})(X)-\overline{y})\cap(-\operatorname*{int}C)=\emptyset.
\]

For several motivations and comments on the minimizing the sum (or difference)
of two (single or set-valued) mappings we refer to the recent works
\cite{Gad1}, \cite{Gad2}. Theorem \ref{T1} allows us to formulate a necessary
optimality condition for $(P).$

\begin{thm}
Let $F_{1},F_{2}:X\rightrightarrows Y,$ be set-valued maps with closed graph,
$(\overline{x},\overline{y}_{1})\in\operatorname*{Gr}F_{1},$ $(\overline
{x},\overline{y}_{2})\in\operatorname*{Gr}F_{2}$ and let $(\overline
{x},\overline{y}_{1}+\overline{y}_{2})$ be a {weak Pareto solution of }$(P).$
Suppose that the function $g:(X\times Y)^{2}\rightarrow X,$ $g(\alpha
,\beta,\gamma,\delta)=\alpha-\gamma$ is metrically subregular at
$(\overline{x},\overline{y}_{1},\overline{x},\overline{y}_{2},0_{X})$ with
respect to $\operatorname*{Gr}F_{1}\times\operatorname*{Gr}F_{2}$.

(i) If either $F_{1}$ is proto-differentiable at $\overline{x}$ relative to
$\overline{y}_{1}$ or $F_{2}$ is proto-differentiable at $\overline{x}$
relative to $\overline{y}_{2}$ then, for every $u\in X,$
\[
\left[  D_{B}F_{1}(\overline{x},\overline{y}_{1})(u)+D_{B}F_{2}(\overline
{x},\overline{y}_{2})(u)\right]  \cap(-\operatorname*{int}C)=\emptyset.
\]

(ii) Let $x\in X$ and $y_{1},y_{2}\in Y$ with $y_{1}+y_{2}\in-C.$ If either
$F_{1}$ is second-order proto-differentiable at $\overline{x}$ relative to
$\overline{y}_{1}$ in the direction $(x,y_{1})$ or $F_{2}$ is second-order
proto-differentiable at $\overline{x}$ relative to $\overline{y}_{2}$ in the
direction $(x,y_{2})$ then for every $u\in X,$%
\[
\left[  D_{B}^{2}F_{1}((\overline{x},\overline{y}_{1})(x,y_{1}))(u)+D_{B}%
^{2}F_{2}((\overline{x},\overline{y}_{2}),(x,y_{2}))(u)\right]  \cap
(-\operatorname*{int}C)=\emptyset.
\]

\end{thm}

\noindent\textbf{Proof.} Since $(\overline{x},\overline{y}_{1}+\overline
{y}_{2})$ is a {weak Pareto solution of }$(P),$ then, following \cite[Chapter
3]{Luc} or \cite[Lemma 3.4]{Dur},
\[
T_{B}((F_{1}+F_{2})(X),\overline{y}_{1}+\overline{y}_{2})\cap
(-\operatorname*{int}C)=\emptyset
\]
and for any $z\in-C$%
\[
T_{B}^{2}((F_{1}+F_{2})(X),\overline{y}_{1}+\overline{y}_{2},z)\cap
(-\operatorname*{int}C)=\emptyset.
\]
But it is easy to see that for any $u\in X,$
\[
D_{B}(F_{1}+F_{2})(\overline{x},\overline{y}_{1}+\overline{y}_{2})(u)\subset
T_{B}((F_{1}+F_{2})(X),\overline{y}_{1}+\overline{y}_{2})
\]
and for any $u\in X,$ $(x_{1},y_{1},y_{2})\in X\times Y\times Y,$%
\[
D_{B}^{2}(F_{1}+F_{2})((\overline{x},\overline{y}_{1}+\overline{y}_{2}%
),(x_{1},y_{1}+y_{2}))(u)\subset T_{B}^{2}((F_{1}+F_{2})(X),\overline{y}%
_{1}+\overline{y}_{2},y_{1}+y_{2}).
\]
Summing up these facts and applying as well Theorem \ref{T1}, with the
identity map instead of $f,$ one gets the conclusions.\hfill$\square$

\bigskip

Moreover, we would like to notice that the nature of our conditions of metric
subregularity in Theorem \ref{T1} are very different from the
semi-differentiability conditions imposed in \cite[Proposition 2]{LiM},
\cite{Li} and \cite[Proposition 2.1]{LeeHuy}. The next example puts the accent
on these differences.

\begin{examp}
Let $F:\mathbb{R}\rightrightarrows\mathbb{R}$ given by
\[
F(x)=\left\{
\begin{array}
[c]{l}%
\emptyset,x=\frac{1}{n},n\in\mathbb{N}\setminus\{0\}\\
x,\text{ otherwise.}%
\end{array}
\right.
\]
where $\mathbb{R}$ and $\mathbb{N}$ denotes the sets of real numbers and
natural numbers, respectively. It is easy to observe that $D_{B}%
F(0,0)(u)=D_{U}F(0,0)(u)=\{u\}$ for any $u\in\mathbb{R}$ but $D_{D}%
F(0,0)(0)=\emptyset$, whence $F$ is proto-differentiable but not
semi-differentiable at $0_{\mathbb{R}}$ relative to $0_{\mathbb{R}}.$ Consider
now $F_{1}=F_{2}:=F.$ The metric subregularity condition in our Theorem
\ref{T1} is nevertheless fulfilled. Moreover, even the metric regularity
around the reference point holds. To see this, take $s>0,$ $p\in(-s,s),$
$\alpha,\gamma\in(-s,s)$. It is enough to show that there is $\mu>0$
(independent of the previous data) s.t.%
\[
\inf\{d((\alpha,\alpha,\gamma,\gamma),(a,a,c,c))\mid a-c=p\}\leq\mu\left\vert
p-(\alpha-\gamma)\right\vert
\]
i.e.
\[
2\inf\{\left\vert \alpha-a\right\vert +\left\vert (\gamma+p)-a\right\vert \mid
a\in\mathbb{R}\}\leq\mu\left\vert (\gamma+p)-\alpha\right\vert
\]
and this is clearly true for any $\mu\geq2.$
\end{examp}

\section{Calculus of derivatives of perturbation maps}

The same technique as in the above section is applied now in order to get
calculus for derivatives of generalized perturbation map. Let $X,Y,Z$ be
Banach spaces, $F,K:X\times Y\rightrightarrows Z$ be set-valued maps$.$ Define
(see \cite{LeeHuy}, \cite{LiM}) $G:X\times Z\rightrightarrows Y$ the
set-valued map given by
\[
G(x,z)=\{y\in Y\mid z\in F(x,y)+K(x,y)\}.
\]
We report now a result concerning a calculus rule for the Bouligand
derivatives of $G.$ Observe first that for $(x,z,y)\in\operatorname*{Gr}G,$
there exist $q_{z}\in F(x,y)$ and $t_{z}\in K(x,y)$ s.t. $z=q_{z}+t_{z}.$ We
use this notation in the sequel.

\begin{thm}
In the above notations, suppose that the graphs of $F$ and $K$ are closed,
$(\overline{x},\overline{z},\overline{y})\in\operatorname*{Gr}G\ $and the
function $i:(X\times Y\times Z)^{2}\rightarrow X\times Y,$
$i(x,y,z,u,v,t)=(x-u,y-v)$ is metrically subregular at $(\overline
{x},\overline{y},q_{\overline{z}},\overline{x},\overline{y},t_{\overline{z}%
},0_{X},0_{Y})$ with respect to $\operatorname*{Gr}F\times\operatorname*{Gr}%
K.$

(i) Then, for every $(u,w)\in X\times Z,$%
\[
\{v\in Y\mid w\in D_{B}F(\overline{x},\overline{y},q_{\overline{z}%
})(u,v)+D_{B}K(\overline{x},\overline{y},t_{\overline{z}})(u,v)\}\subset
D_{B}G(\overline{x},\overline{z},\overline{y})(u,w),
\]

(ii) Let $(u,v,t)\in X\times Y\times Z$ with $t=t_{1}+t_{2}$ $(t_{1},t_{2}\in
Z)$ and suppose that either $F\ $is second-order proto-differentiable at
$(\overline{x},\overline{y})$ relative to $q_{\overline{z}}$ in the direction
$(u,v,t_{1})$ or $K\ $is second-order proto-differentiable at $(\overline
{x},\overline{y})$ relative to $t_{\overline{z}}$ in the direction
$(u,v,t_{2}).$ Then for any $(\alpha,\beta)\in X\times Y$ one has%
\begin{align*}
\gamma_{1}  &  \in D_{B}^{2}F((\overline{x},\overline{y},q_{\overline{z}%
}),(u,v,t_{1}))(\alpha,\beta),\gamma_{2}\in D_{B}^{2}F((\overline{x}%
,\overline{y},q_{\overline{z}}),(u,v,t_{1}))(\alpha,\beta)\\
&  \Rightarrow\beta\in\operatorname*{Gr}D_{B}^{2}G(\overline{x},\overline
{z},\overline{y})(u,t,v)(\alpha,\gamma_{1}+\gamma_{2}).
\end{align*}

\end{thm}

\noindent\textbf{Proof.} (i) Take
\[
\varphi:(X\times Y\times Z)^{2}\rightarrow X\times Z\times Y,\quad
\varphi(x,y,z,u,v,t)=(x,z+t,y)
\]
and%
\[
\psi:(X\times Y\times Z)^{2}\rightarrow X\times Y,\quad\psi
(x,y,z,u,v,t)=(x-u,y-v).
\]
Then
\begin{align*}
&  \varphi((\operatorname*{Gr}F\times\operatorname*{Gr}K)\cap\psi^{-1}%
(0_{X},0_{Y}))\\
&  =\varphi(\{(x,y,z,u,v,t)\mid(x,y,z)\in\operatorname*{Gr}F,(u,v,t)\in
\operatorname*{Gr}K,x=u,y=v\})\\
&  =\varphi(\{(x,y,z,x,y,t)\mid(x,y,z)\in\operatorname*{Gr}F,(x,y,t)\in
\operatorname*{Gr}K\})\\
&  =\{(x,z+t,y)\mid z\in F(x,y),t\in K(x,y)\}\\
&  =\{(x,w,y)\mid w\in F(x,y)+K(x,y)\}=\operatorname*{Gr}G.
\end{align*}
Therefore,
\begin{align*}
\operatorname*{Gr}D_{B}G(\overline{x},\overline{z},\overline{y}) &
=T_{B}(\operatorname*{Gr}G,(\overline{x},\overline{z},\overline{y}))\\
&  =T_{B}(\varphi((\operatorname*{Gr}F\times\operatorname*{Gr}K)\cap\psi
^{-1}(0_{X},0_{Y})),(\overline{x},q_{\overline{z}}+t_{\overline{z}}%
,\overline{y}))\\
&  =T_{B}(\varphi((\operatorname*{Gr}F\times\operatorname*{Gr}K)\cap\psi
^{-1}(0_{X},0_{Y})),\varphi(\overline{x},\overline{y},q_{\overline{z}%
},\overline{x},\overline{y},t_{\overline{z}}))\\
&  \supset\operatorname*{cl}\left(  \varphi(T_{B}((\operatorname*{Gr}%
F\times\operatorname*{Gr}K)\cap\psi^{-1}(0_{X},0_{Y}),(\overline{x}%
,\overline{y},q_{\overline{z}},\overline{x},\overline{y},t_{\overline{z}%
})))\right)  .
\end{align*}
Using the metric subregularity assumption, the function $j:(X\times Y\times
Z)^{2}\times(X\times Y)\rightarrow X\times Y,$
$j(x,y,z,u,v,t,p,q)=(x-u-p,y-v-q)$ is metrically subregular at $(\overline
{x},\overline{y},q_{\overline{z}},\overline{x},\overline{y},t_{\overline{z}%
},0_{X},0_{Y},0_{X},0_{Y})$ with respect to $(\operatorname*{Gr}%
F\times\operatorname*{Gr}K)\times\{(0_{X},0_{Y})\}\ $whence applying Theorem
\ref{th_teor} one has:
\begin{align*}
&  T_{B}((\operatorname*{Gr}F\times\operatorname*{Gr}K)\cap\psi^{-1}%
(0_{X},0_{Y}),(\overline{x},\overline{y},q_{\overline{z}},\overline
{x},\overline{y},t_{\overline{z}}))\\
&  \supset T_{B}(\operatorname*{Gr}F\times\operatorname*{Gr}K,(\overline
{x},\overline{y},q_{\overline{z}},\overline{x},\overline{y},t_{\overline{z}%
}))\cap\nabla\psi(\overline{x},\overline{y},q_{\overline{z}},\overline
{x},\overline{y},t_{\overline{z}})^{-1}(0_{X},0_{Y})).
\end{align*}
We employ now the proto-differentiability assumption and, accordingly, one
has
\begin{align*}
\operatorname*{Gr}D_{B}G(\overline{x},\overline{z},\overline{y}) &
\supset\varphi(\left(  T_{B}(\operatorname*{Gr}F,(\overline{x},\overline
{y},q_{\overline{z}}))\times T_{B}(\operatorname*{Gr}K,(x,y,t_{\overline{z}%
}))\right)  \cap\nabla\psi(\overline{x},\overline{y},q_{\overline{z}%
},\overline{x},\overline{y},t_{\overline{z}})^{-1}(0_{X},0_{Y}))\\
&  =\varphi(\{(u,v,p,r,s,j)\mid p\in D_{B}F(\overline{x},\overline
{y},q_{\overline{z}})(u,v),j\in D_{B}K(\overline{x},\overline{y}%
,t_{\overline{z}})(r,s),u=r,v=s\})\\
&  =\{(u,p+j,v)\mid p\in D_{B}F(\overline{x},\overline{y},q_{\overline{z}%
})(u,v),j\in D_{B}K(\overline{x},\overline{y},t_{\overline{z}})(u,v)\}.
\end{align*}
We conclude that
\[
\{v\in Y\mid w\in D_{B}F(\overline{x},\overline{y},q_{\overline{z}%
})(u,v)+D_{B}K(\overline{x},\overline{y},t_{\overline{z}})(u,v)\}\subset
D_{B}G(\overline{x},\overline{z},\overline{y})(u,w),
\]
hence the first conclusion.

(ii) For the second part, we firstly write
\begin{align*}
&  \operatorname*{Gr}D_{B}^{2}G(\overline{x},\overline{z},\overline
{y})(u,t,v)=\\
&  =T_{B}^{2}(\operatorname*{Gr}G,(\overline{x},\overline{z},\overline
{y})(u,t,v))\\
&  =T_{B}^{2}(\varphi((\operatorname*{Gr}F\times\operatorname*{Gr}K)\cap
\psi^{-1}(0_{X},0_{Y})),\varphi(\overline{x},\overline{y},q_{\overline{z}%
},\overline{x},\overline{y},t_{\overline{z}}),\varphi(u,v,t_{1},u,v,t_{2}))\\
&  \supset\varphi\left(  T_{B}^{2}(\operatorname*{Gr}F\times\operatorname*{Gr}%
K)\cap\psi^{-1}(0_{X},0_{Y}),(\overline{x},\overline{y},q_{\overline{z}%
},\overline{x},\overline{y},t_{\overline{z}}),(u,v,t_{1},u,v,t_{2}))\right)  .
\end{align*}

By similar arguments,
\begin{align*}
&  T_{B}^{2}(\operatorname*{Gr}F\times\operatorname*{Gr}K)\cap\psi^{-1}%
(0_{X},0_{Y}),(\overline{x},\overline{y},q_{\overline{z}},\overline
{x},\overline{y},t_{\overline{z}}),(u,v,t_{1},u,v,t_{2}))\\
&  \supset T_{B}^{2}(\operatorname*{Gr}F\times\operatorname*{Gr}%
K,(\overline{x},\overline{y},q_{\overline{z}},\overline{x},\overline
{y},t_{\overline{z}}),(u,v,t_{1},u,v,t_{2}))\\
&  \cap\nabla\psi(\overline{x},\overline{y},q_{\overline{z}},\overline
{x},\overline{y},t_{\overline{z}})^{-1}(0_{X},0_{Y}).
\end{align*}

The second-order proto-differentiability condition ensures
\begin{align*}
&  \operatorname*{Gr}D_{B}^{2}G(\overline{x},\overline{z},\overline
{y})(u,t,v)\\
&  \supset\varphi\left(  T_{B}^{2}(\operatorname*{Gr}F\times\operatorname*{Gr}%
K,(\overline{x},\overline{y},q_{\overline{z}},\overline{x},\overline
{y},t_{\overline{z}}),(u,v,t_{1},u,v,t_{2}))\cap\nabla\psi(\overline
{x},\overline{y},q_{\overline{z}},\overline{x},\overline{y},t_{\overline{z}%
})^{-1}(0_{X},0_{Y})\right)  \\
&  =\varphi\left(  T_{B}^{2}(\operatorname*{Gr}F,(\overline{x},\overline
{y},q_{\overline{z}}),(u,v,t_{1}))\times T_{B}^{2}(\operatorname*{Gr}%
K,(\overline{x},\overline{y},t_{\overline{z}}),(u,v,t_{2}\right)  )\\
&  \hspace{1.9cm}\cap\nabla\psi(\overline{x},\overline{y},q_{\overline{z}%
},\overline{x},\overline{y},t_{\overline{z}})^{-1}(0_{X},0_{Y}))\\
&  =\varphi(\{(a,b,c,d,e,f)\in(X\times Y\times Z)^{2}\mid(a,b,c)\in
\operatorname*{Gr}D_{B}^{2}F((\overline{x},\overline{y},q_{\overline{z}%
}),(u,v,t_{1})),\\
&  \hspace{1.9cm}(d,e,f)\in\operatorname*{Gr}D_{B}^{2}K((\overline
{x},\overline{y},t_{\overline{z}}),(u,v,t_{2})),a=d,b=e\}\\
&  =\{(a,c+f,b)\mid(a,b,c)\in\operatorname*{Gr}D_{B}^{2}F((\overline
{x},\overline{y},q_{\overline{z}}),(u,v,t_{1})),\\
&  \hspace{1.9cm}(a,b,f)\in\operatorname*{Gr}D_{B}^{2}K((\overline
{x},\overline{y},t_{\overline{z}}),(u,v,t_{2}))\}.
\end{align*}

Consequently, the second part of the conclusion easily follows.\hfill$\square$

\medskip

\begin{rmk}
On the basis of Corollary \ref{cor_mr}, a sufficient condition for the metric
(sub)regularity assumption in Theorem \ref{T1} (on Asplund spaces) reads as
follows: there exist $c>0,$ $r>0$ such that for every $(x_{1},y_{1}%
,z_{1},x_{2},y_{2},z_{2})\in(\operatorname*{Gr}F\times\operatorname*{Gr}%
K)\cap\lbrack B(\overline{x},r)\times B(\overline{y},r)\times B(q_{\overline
{z}},r)\times B(\overline{x},r)\times B(\overline{y},r)\times B(t_{\overline
{z}},r)]$ and every $(x^{\ast},y^{\ast})\in X^{\ast}\times Y^{\ast},$
$(u_{1}^{\ast},v_{1}^{\ast},w_{1}^{\ast},u_{2}^{\ast},v_{2}^{\ast},w_{2}%
^{\ast})\in\widehat{N}(\operatorname*{Gr}F,(x_{1},y_{1},z_{1}))\times
\widehat{N}(\operatorname*{Gr}K,(x_{2},y_{2},z_{2}))+(x^{\ast},y^{\ast
},0_{Z^{\ast}},-x^{\ast},-y^{\ast},0_{Z^{\ast}})$%
\[
c\left\Vert (x^{\ast},y^{\ast})\right\Vert \leq\left\Vert (u_{1}^{\ast}%
,v_{1}^{\ast},w_{1}^{\ast},u_{2}^{\ast},v_{2}^{\ast},w_{2}^{\ast})\right\Vert
.
\]

\end{rmk}

\bigskip

We consider now a generalization of the perturbation map from the
Rockafellar's paper \cite{Roc} (see \cite[Theorem 5.4]{Roc}). Let, as above,
$X,Y,Z$ be general Banach spaces, $f:X\times Y\rightarrow Z$ be a continuously
differentiable single-valued map, $D\subset Y,E\subset Z$ be closed subsets
and define the set-valued $F:X\rightrightarrows Y$ defined by:%
\[
F(x)=\{y\in D\mid f(x,y)\in E\}.
\]
The calculus of the derivatives of this model set-valued map can also be
captured by the approach we use in this paper.

\begin{thm}
In the above notations, suppose that $D,E$ are closed, $(\overline
{x},\overline{y})\in X\times Y$ with $\overline{x}\in D,f(\overline
{x},\overline{y})\in E$ and $g:X\times Y\times Z\longrightarrow Z,$
$g(x,y,z)=f(x,y)-z$ is metrically subregular at $(\overline{x},\overline
{y},f(\overline{x},\overline{y}),0_{Z})$ with respect to $(X\times D)\times
E.$

(i) If either $D$ is derivable at $\overline{y}$ or $E$ is derivable at
$f(\overline{x},\overline{y})$, then for every $u\in X,$
\[
D_{B}F(\overline{x},\overline{y})(u)=\{v\in T_{B}(D,\overline{y})\mid\nabla
f(\overline{x},\overline{y})(u,v)\in T_{B}(E,f(\overline{x},\overline{y}))\}.
\]

(ii) Suppose that $f$ is twice continuously differentiable. If either $D$ is
second-order derivable at $\overline{y}$ or $E$ is second-order derivable at
$f(\overline{x},\overline{y}),$ then for every $x_{1}\in X$ and $y_{1}\in Y$%
\begin{align*}
\{v  &  \in T_{B}^{2}(D,\overline{y},y_{1})\mid\nabla f(\overline{x}%
,\overline{y})(u,v)\in T_{B}^{2}(E,f(\overline{x},\overline{y}),\nabla
f(\overline{x},\overline{y})(x_{1},y_{1}))\\
&  -2^{-1}\nabla^{2}f(\overline{x},\overline{y})((x_{1}y_{1}),(x_{1}%
,y_{1}))\}\\
&  =D_{B}^{2}F(\overline{x},\overline{y})(x_{1},y_{1})
\end{align*}

\end{thm}

\noindent\textbf{Proof.} (i) Obviously, $\operatorname*{Gr}F=(X\times D)\cap
f^{-1}(E).$ Then, using Proposition \ref{prop_con} (ii),
\begin{align*}
T_{B}(\operatorname*{Gr}F,(\overline{x},\overline{y}))  &  =T_{B}((X\times
D)\cap f^{-1}(E),(\overline{x},\overline{y}))\\
&  \subset T_{B}(X\times D,(\overline{x},\overline{y}))\cap\nabla
f(\overline{x},\overline{y})^{-1}(T_{B}(E,f(\overline{x},\overline{y}))),
\end{align*}
Clearly, $T_{B}(X\times D,(\overline{x},\overline{y}))=X\times T_{B}%
(D,\overline{y})$ (see Proposition \ref{prop_con} (i)), so%
\[
T_{B}(\operatorname*{Gr}F,(\overline{x},\overline{y}))\subset\left(  X\times
T_{B}(D,\overline{y})\right)  \cap\nabla f(\overline{x},\overline{y}%
)^{-1}(T_{B}(E,f(\overline{x},\overline{y}))).
\]
In order to prove the converse inclusion observe that in our hypotheses one
can apply Theorem \ref{th_teor} for $X\times D,E$ and $f.$ Therefore, using as
well the derivability assumption on the sets,
\[
T_{B}((X\times D)\cap f^{-1}(E),(x,y))\supset\left(  X\times T_{B}%
(D,y)\right)  \cap\nabla f(x,y)^{-1}(T_{B}(E,f(x,y))).
\]
The conclusion follows.

(ii) This time, the second part of Theorem \ref{th_teor} yields the
conclusion:%
\begin{align*}
&  T_{B}^{2}(X\times D,(\overline{x},\overline{y}),(x_{1},y_{1}))\cap\nabla
f(\overline{x},\overline{y})^{-1}(T_{B}^{2}(E,f(\overline{x},\overline
{y}),\nabla f(\overline{x},\overline{y})(x_{1},y_{1}))\\
&  -2^{-1}\nabla^{2}f(\overline{x},\overline{y})((x_{1}y_{1}),(x_{1}%
,y_{1})))\\
&  \subset T_{B}^{2}((X\times D)\cap f^{-1}(E),(\overline{x},\overline
{y}),(x_{1},y_{1})).
\end{align*}

The fact that $T_{B}^{2}(X\times D,(\overline{x},\overline{y}),(x_{1}%
,y_{1}))=X\times T_{B}^{2}(D,\overline{y},y_{1})$ ends the proof of this part.

The other inclusion is a simple application of second-order Taylor
formula.$\hfill\square$

\bigskip

\begin{rmk}
As we have proceeded before, we give now (on Asplund spaces) the condition for
the metric (sub)regularity, according to Corollary \ref{cor_mr}: there exist
$c>0,$ $r>0$ such that for every $(x,y,z)\in\lbrack X\times D\times
E]\cap\lbrack B(\overline{x},r)\times B(\overline{y},r)\times B(f(\overline
{x},\overline{y}),r)]$ and every $z^{\ast}\in Z^{\ast},$ $(u^{\ast},v^{\ast
},w^{\ast})\in X\times\widehat{N}(D,y)\times\widehat{N}(E,z)+(\nabla^{\ast
}f(x,y)(z^{\ast}),-z^{\ast})$%
\[
c\left\Vert z^{\ast}\right\Vert \leq\left\Vert (u^{\ast},v^{\ast},w^{\ast
})\right\Vert .
\]

\end{rmk}

\section{Concluding remarks}

This paper applies a method based on metric subregularity assumptions in order
to compute the first and second order derivatives of several set-valued maps.
This approach is natural in infinite dimensional setting since allows us to
eliminate any compactness requirement. Moreover, we have described, on Asplund
spaces, sufficient conditions for our hypotheses in terms of Fr\'{e}chet
normal cones and we have stressed that similar conditions could be imposed as
well on other appropriate classes of Banach spaces and corresponding
well-behaved normals. In the case of set-valued mappings with closed convex
graphs these limitation can be dropped. We have mainly concentrated on the sum
of set-valued mappings and on perturbation maps. Nevertheless, we would like
to emphasize that this method could be applied in further investigations for
other models as well and the inclusion obtained in such a way could be useful
in the study of other vector optimization programs, different from the one
considered here.

\end{document}